\documentclass[amsmath, amsthm, amssymb, floatfix, reprint, twocolumn, notitlepage, nofootinbib, 10pt]{revtex4-1}

\newtheorem{theorem}{Theorem}
\newtheorem{lemma}{Lemma}
\newtheorem{problem}{Problem}
\newtheorem{proposition}{Proposition}

\newtheorem{definition}{Definition}

\usepackage[utf8]{inputenc}
\usepackage{microtype,bm,bbm,graphicx,booktabs,times}
\usepackage[usenames,dvipsnames]{xcolor}
\usepackage{setspace}
\usepackage{algorithm}
\usepackage{algorithmic}
\selectfont{}

\usepackage{hyperref}
\usepackage{braket}
\usepackage{subfiles} 
\usepackage{upgreek}
\usepackage{color}
\hypersetup{colorlinks,allcolors=black,linktoc=all}
\usepackage[capitalize]{cleveref}
\crefformat{section}{Sec.~#2#1#3}

\newcommand{\norm}[1]{\Vert{#1}\Vert}

\begin{document}
\title{Gradient descent globally solves average-case non-resonant physical design problems}
\author{Rahul Trivedi}
\email{rahul.trivedi@mpq.mpg.de}
\address{Max-Planck-Institut für Quantenoptik, Hans-Kopfermann-Str.~1, 85748 Garching, Germany.}

\date{\today}

\begin{abstract}
Optimization problems occurring in a wide variety of physical design problems, including but not limited to optical engineering, quantum control, structural engineering, involve minimization of a simple cost function of the state of the system (e.g.~the optical fields, the quantum state) while being constrained by the physics of the system. The physics constraints often makes such problems non-convex and thus only locally solvable, leaving open the question of finding the globally optimal design. In this paper, I consider design problems whose physics is described by bi-affine equality constraints, and show that under assumptions on the stability of these constraints and the physical system being non-resonant, gradient descent globally solves a typical physical design problem. The key technical contributions of this paper are (i) outline a criteria that ensure the convergence of gradient descent to an approximate global optima in the limit of large problem sizes, and (ii) use random matrix theory to outline ensembles of physically motivated problems which, on an average, satisfy this convergence criteria. 
\end{abstract}

\maketitle
\section{Introduction}
Optimization-based design methodologies have been successfully adapted and applied to solving a wide variety of physical design problems, including but not limited to photonics design \cite{molesky2018inverse, su2018inverse, piggott2019inverse, su2018fully, piggott2017fabrication, sapra2019inverse}, device-level designs for quantum technologies \cite{roslund2009gradient, lucarelli2018quantum, werschnik2007quantum, li2017hybrid, zhu1998rapid}, and structural engineering \cite{rozvany2009critical, wang2003level, wang2003level,  hassani2012homogenization}. The physics of these system are mathematically captured by bi-affine equality constraints, i.e.~constraints of the form $h(x, \theta) = 0$ where $h$ is affine in $x$ for fixed $\theta$ and affine in $\theta$ for fixed $x$, relating the system state ($x$) to the design parameters ($\theta$). The design problem is captured by a cost function of the system state $x$, which is to be minimized with respect to the parameters $\theta$ while enforcing these constraints. In practice, the bi-affine equality constraint allow for an efficient evaluation of the gradient of the cost function with respect to the design parameters \cite{giles2000introduction,bakr2004adjoint,  park1996design, dong2004design, soliman2004adjoint}, and thus local optimization algorithms such as gradient descent \cite{ruder2016overview} of quasi-Newton methods \cite{gill1972quasi} can be used to solve them locally. However, the bi-affine constraints also makes the optimization problem non-convex, and hence hard to solve globally. While the locally optimized result suffices for many applications, it leaves open the question of how much more improvement in the device performance can be gained by searching for the global optima.

One prominent approach to answering this question has been to provide lower bounds on the cost function that captures the device performance and is being optimized --- while the optimization problem itself is non-convex and hard to solve globally, such bounds can often be computed efficiently by solving a convex problem. Several approaches to setting up this convex problem and calculating these lower bounds using physically motivated convex relaxations \cite{miller2016fundamental, shim2019fundamental, venkataram2020fundamental, molesky2020fundamental, molesky2019t} or by an application of Lagrange duality \cite{angeris2019computational,kuang2020maximal, kuang2020computational, prl_qcontrol} have been recently pursued. In many problems of interest, these lower bounds, computed numerically, are reasonably close to the locally optimized results and thus indicate that the design is near globally-optimal \cite{zhang2020optimal, kern2020tight}. However, there is no general theoretical guarantee for the lower bounds obtained by these methods to be tight (i.e.~close to the global optima), and consequently they do not always allow for a quantitative assessment of the optimality of the locally optimized design.

In this paper, I pursue a different approach for assessing the optimality of local optimization algorithms --- I theoretically analyze the convergence of the gradient descent algorithm when applied to a physical design problem. The results indicate that gradient descent efficiently solves a typical physical design problem i.e.~it comes within a specified accuracy of the global optima in a number of steps that only scale polynomially with the problem size. There are two parts to the main result --- first, I provide a set of conditions on the physical design problem which guarantee this convergence result. Next, I analyze physically motivated distributions of bi-affine design problems and show that a  member of this ensemble, on average, satisfies these conditions and thus is expected to be efficiently solvable by gradient descent. While this work is, to the best of my knowledge, the first rigorous study of convergence of gradient descent for physical design problems, it shares techniques with similar analysis of local optimization algorithms when applied to non-convex problems arising in training of deep neural networks \cite{allen2019convergence, allen2019learning, allen2019convergence}, neural-tangent kernels \cite{jacot2021neural} and dynamical models of time-series data \cite{hardt2018gradient}.

\section{Notation}
For $v \in \mathbb{R}^n$, I will denote by $\norm{v}_p$ its $l^p$ norm, and for simplicity use $\norm{v}$ for its $l^2$ norm. For a matrix $M\in \mathbb{R}^{n\times m}$, I will denote by $\norm{M}_p$ the induced $l^p$ norm i.e.~$\norm{M}_p = \sup_{v\in \mathbb{R}^m \setminus\{0\}} \norm{Mv}_p / \norm{v}_p$.  I will denote by $I_n$ the $n\times n$ identity matrix. I will denote by $\norm{M}_\textnormal{max}$ the maximum magnitude element of $M$ i.e.~$\norm{M}_\textnormal{max} = \max_{i \in [n], j \in [m]} |M_{i, j}|$. I will often write the vector, $0^n$ (i.e.~vector with all elements 0), simply as 0, and the dimensionality of the vector will be evident from the context.

For two vectors $u, v \in \mathbb{R}^n$, $u \odot v \in \mathbb{R}^n$ is their elementwise product defined by $(u \odot v)_i := u_i v_i \ \forall \ i \in [n]$.  Given a function $f: \text{dom}(f) \subset \mathbb{R} \to \mathbb{R}$ and a vector $v \in \text{dom}(f)^n$, $f(v)$ will denote a vector obtained on applying $f$ entry-wise to $f$.

I will use the computer science notation for asymptotic behaviour of sequences. Given a sequence $\{a_n \in \mathbb{R}^+ : n \in \mathbb{N}\}$. $a_n \leq O(f(n))$ for some $f: \mathbb{N} \to \mathbb{R}^+$ if $\exists c > 0$ such that $a_n \leq c f(n)$ as $n \to \infty$. $a_n < o(f(n))$ for some $f: \mathbb{N} \to \mathbb{R}^+$ if $\forall c > 0$, $a_n < c f(n)$ as $n \to \infty$. In particular, $n^{o(1)}$ will be used to denote a function whose growth is slower than  $n^\alpha$ for any $\alpha > 0$. $a_n \geq \Omega(f(n))$ for some $f: \mathbb{N} \to \mathbb{R}^+$ if $\exists c>0$ such that $a_n \geq c f(n)$ as $n \to \infty$. $a_n = \Theta(f(n))$ if $a_n \leq O(f(n))$ and $a_n \geq \Omega(f(n))$.

\section{Summary of results}

I first introduce a definition of a physical design problem, along with some terminology that I use throughout this paper. I introduce an abstract definition of a \emph{physical system}, and a \emph{physical design problem}.
\begin{definition} A \emph{physical system} with state size $n$ and parameter size $m$ is a map $\varphi:\textnormal{dom}(\varphi) \to \mathbb{R}^n$ specified by a tuple $(A, B, b)$ where $A \in \mathbb{R}^{n \times n}$, $B \in \mathbb{R}^{n\times m}$, $b \in \mathbb{R}^n$ and
\[
\textnormal{dom}(\varphi) = \{ \theta \in \mathbb{R}^m | A + \textnormal{diag}(B\theta) \textnormal{ is invertible}\},
\]
and $\forall \theta \in \mathbb{R}^m$,
\[
\varphi(\theta) = (A + \textnormal{diag}(B\theta))^{-1}b.
\]
$A$ is the physics matrix, $B$ is the selection matrix and $b$ is the source vector corresponding to the physical system.
\end{definition}
\begin{definition}Given a physical system $\varphi \equiv (A, B, b)$ of state size $n$, and a vector $c \in \mathbb{R}^n$, the \emph{adjoint system} is a map $\textnormal{ad}[\varphi]: \textnormal{dom}(\varphi) \times \mathbb{R}^n \to \mathbb{R}^n$ such that $\forall \theta \in \mathbb{R}^m$,
\[
\textnormal{ad}[\varphi](\theta, c) = (A + \textnormal{diag}(B \theta))^{-\textnormal{T}} c.
\]
\end{definition}
 
\begin{definition}
A physical design problem of state size $n$ and parameter size $m$ is specified by a tuple $(f, c, \varphi)$ where $f: \mathbb{R} \to \mathbb{R}$ is the cost function, $c \in \mathbb{R}^n$ is the overlap vector, $\varphi$ is a physical system of state size $n$ and parameter size $m$ and it corresponds to solving the following constrained optimization problem:
\begin{equation*}\label{eq:global_problem}
\begin{aligned}
& \underset{\substack{\theta \in \textnormal{dom}(\varphi)}}{\textnormal{minimize}} & & f(c^\textnormal{T}\varphi(\theta)),\\
\end{aligned}
\end{equation*}
\end{definition}

To solve this problem using local optimization algorithm, it is essential to be able to efficiently compute the gradient of the cost function. It is well known that this can be done using the the maps corresponding to the physical system and the adjoint system \cite{giles2000introduction,bakr2004adjoint,  park1996design, dong2004design, soliman2004adjoint}. This is made explicit in the following lemma, which can be straightforwardly proved using the chain rule.
\begin{lemma}
\label{lemma:adj_method}
For a physical design problem $(f, c, \varphi)$, then for the map $f(c^\textnormal{T}\varphi(\cdot)):\textnormal{dom}(\varphi) \to \mathbb{R}$, the gradient at $\theta \in \textnormal{dom}(\varphi)$ is given by
\[
\nabla_\theta f(c^\textnormal{T}\varphi(\theta)) = -f'(c^\textnormal{T}\varphi(\theta)) B_\varphi^\textnormal{T} (\varphi(\theta) \odot \textnormal{ad}[\varphi](\theta, c)).
\]
\end{lemma}

\begin{algorithm}[H] 
\caption{Gradient descent for solving a physical design problem}
\label{alg:gd}
\begin{algorithmic}
\renewcommand{\algorithmicrequire}{\textbf{Input:}}
 \renewcommand{\algorithmicensure}{\textbf{Output:}}
 \REQUIRE A physical design problem $(f, c, \varphi)$, an initial set of design parameters $\theta_0 \in \mathbb{R}^m$, a gradient descent step size $\eta \in (0, \infty)$ and number of gradient descent steps $T \in \mathbb{N}$.
 \ENSURE  The optimized design parameters $\theta^* \in \mathbb{R}^m$.
 \\ \textit{Initialisation} :
  \STATE If $\theta_0 \notin \textnormal{dom}(\varphi)$, then declare FAIL, else set $x_0 = \varphi(\theta_0) = (A_\varphi + \textnormal{diag}(B_\varphi \theta))^{-1} b_\varphi$. 
 \\ \textit{LOOP Process}
  \FOR {$t = 1$ to $T - 1$}
  \STATE Compute the adjoint state $a_{t - 1} = \varphi^\textnormal{adj}_c (\theta_{t - 1})$.
  \STATE Compute the gradient $g_{t - 1} = -f'(c^\textnormal{T}x_{t - 1}) \big(\varphi^\textnormal{adj}_c (\theta_{t - 1}) \odot \varphi(\theta_{t-1})\big)$
  \STATE Perform the gradient descent step $\theta_t:= \theta_{t - 1} - \eta g_{t - 1}$.
  \IF {($\theta_t \notin \textnormal{dom}(\varphi)$)}
  \STATE Declare FAIL.
    \ENDIF
  \ENDFOR
 \RETURN $\theta_T$
\end{algorithmic}
\end{algorithm}
Throughout this paper, I will be interested in performing an analysis of gradient descent (explicitly described in algorithm 1) not for one specific problem instance, but for problem instance with large state sizes. More formally, I will consider a family of physical design problems with larger and larger state sizes and assess how gradient descent performs on a specific instance of this family.
\begin{definition} A \emph{family of physical design problems} is a sequence $\{(f, c_n, \varphi_n): n \in \mathbb{N}\}$ of physical design problems, where $(f, c_n, \varphi_n)$ is a physical design problem of state size $n$ and the cost function $f:\mathbb{R} \to \mathbb{R}$ is independent of $n$.
\end{definition}

I next introduce a set of conditions on a given family of physical design problems, which will be central to the analysis of the convergence of gradient descent. These conditions are expressed in terms of the scalings of the spectrum of the physics matrix with the state size of the physical system, as well as on the norms of the source and overlap vectors. The definition I provide below makes concrete the physical expectation that the norm of the physics matrix of systems encountered in practice only grows polynomially with the state size of the system, and typically the source and overlap vectors only affect a constant number of state variables and hence have norms upper bounded by a constant. I also assume a condition on the initial gradient of the cost function --- since in practice the initial design is picked randomly and is thus not locally optimal, we assume that the initial gradient of the cost function is large (i.e.~lower bounded by a function that grows polynomially with the state size). 
\begin{definition}
\label{def:asym_conv_cond}
A family of physical design problems $\{(f, c_n, \varphi_n \equiv (A_n, B_n, b_n)): n \in \mathbb{N}\}$ is said to satisfy the $(\alpha, \gamma)$-asymptotic convergence conditions if
\begin{itemize}
\item[(a)] $f$ is $L-$Lipschitz smooth\footnote{A differentiable function $f:\mathbb{R}^D \to \mathbb{R}$ is $L-$Lipschitz smooth if $\norm{\nabla f(x) - \nabla f(y)} \leq L\norm{x - y}$}, $\mu-$strongly convex \footnote{A function $f:\mathbb{R}^D \to \mathbb{R}$ is $\mu-$strongly convex if $f(y) \geq f(x) + \nabla f(x)^\textnormal{T}(y - x) + \sigma \norm{x - y}^2 / 2$.} and is bounded from below i.e.~$f^* = \min_{x\in \mathbb{R}}f(x) > -\infty$.
\item[(b)] $\norm{A_{n}^{-1}}_2 \leq O(n^\gamma)$,
\item[(c)] $ \norm{B_{n}}_\infty, \norm{B_{n}}_1, \norm{B_{n}}_2 \leq O(1)$,
\item[(d)] $\norm{b_{n}}_\infty, \norm{b_{n}}_1 \leq O(1)$,
\item[(e)] $\norm{c_n}_\infty, \norm{c_n}_1 \leq O(1)$,
\item[(f)] $\norm{B_{n}^\textnormal{T}(A_{n}^{-1}b_{n} \odot A_{n}^{-\textnormal{T}}c_n)}_2^2 \geq \Omega(n^{\alpha + \gamma})$ and
\item[(g)] $|c_{n}^\textnormal{T}A_{n}^{-1}b_{n}| \leq O(1)$.
\end{itemize}
\end{definition}

Typical physical design problems are expected to be non-resonant i.e.~the state vector is expected to not become very large during the optimization trajectory. It is empirically observed in many settings that local optimization algorithms perform very well when applied to such non-resonant design problems, while resonant design problems are usually much harder to solve. This consideration has to be accounted for in the analysis of gradient descent --- the next definition that I provide makes the intuition behind a non-resonant design mathematically precise, and will be another important ingredient in the analysis of the optimality of gradient descent.

\begin{definition} [Non-resonant design parameters] Given a family of physical design problems $\mathcal{F} = \{(f, c_n, \varphi_n ): n \in \mathbb{N}\}$, a sequence of design parameters $\{\theta_n \in \textnormal{dom}(\varphi_n) : n \in \mathbb{N}\}$ is said to be non-resonant with respect to the family $\mathcal{F}$ if $\norm{\varphi_n(\theta_n)}_\infty \leq O(n^{o(1)})$ and $\norm{\textnormal{adj}[\varphi_n](\theta_n, c_n)}_\infty \leq O(n^{o(1)})$.
\end{definition}
I now present the first result of this paper --- a family of physical design problems is efficiently solvable by gradient descent, under the assumption that the design parameters generated during the algorithm are non-resonant, if it satisfies the asymptotic convergence conditions.

\begin{theorem}
\label{thm:conv_grad_des_non_res}
Let $\mathcal{F} := \{(f, c_n, \varphi_n) : n \in \mathbb{N}\}$ be a family of physical design problems that satisfies the $\alpha, \gamma-$asymptotic convergence conditions (definition 5) with $\alpha > 1/4, \gamma < 3\alpha$, and for $n \in \mathbb{N}$, let $\Theta_n = \{\theta_{n}^1, \theta_{n}^2 \dots \theta_{n}^{T_n} \}$ be the design parameters generated by gradient descent (algorithm 1) under the assumption that it does not fail when applied on $(f, c_n, \varphi_n)$ with step size $\eta_n \leq O(n^{-3\gamma + \alpha - 1})$ for $T_n$ steps starting with initial design of $\theta_n^0 = 0$. If all the sequences $\{\theta_{n} \in \Theta_n: n \in \mathbb{N}\}$ are non-resonant with respect to $\mathcal{F}$, then $f(c_n^\textnormal{T}\varphi_n(\theta_{n}^{T_n})) - f^* \leq \varepsilon$ for $T_n$ chosen such that $T_n = \Theta(n^{1 - 2(\alpha - \gamma)}\log(\varepsilon^{-1}))$.
\end{theorem}

In the next two propositions, I provide families of design problems which are constructed by choosing the physics matrices randomly from a distribution. By analyzing an average-case problem picked from these distributions, I show that it satisfies the asymptotic convergence conditions. The first family of problem is one in which the inverse of the physics matrix is a matrix from the random gaussian ensemble, and this example is inspired from wave design problems \cite{molesky2018inverse, su2018inverse, piggott2019inverse, su2018fully}, where often the initial design is a random scattering media whose properties are known to be captured by random matrices \cite{kogan1995random}. 

\begin{problem}
\label{prob:gaussian_random}
Given a cost function $f:\mathbb{R}\to \mathbb{R}$ being an $L-$Lipschitz smooth and $\mu-$strongly convex function that is bounded from below, define a family of physical design problems $\mathcal{F}:= \{(f, c_n, \varphi_n \equiv (A_n: = G_n^{-1},B_n:= I_n, b_n)) : n \in \mathbb{N}\}$ where
\begin{itemize}
\item $\norm{b_n}_\infty, \norm{b_n}_1, \norm{c_n}_\infty, \norm{c_n}_1 = \Theta(1)$ and,
\item $G_n \in \mathbb{R}^{n\times n}$ is a matrix where each entry is independently drawn from the standard normal distribution.
\end{itemize}
\end{problem}

\begin{proposition}  For the problem of state size $n$, $(f, c_n, \varphi_n \equiv (A_n , B_n, b_n))$, picked from the family of problems defined in problem \ref{prob:gaussian_random}, it is true that
\begin{enumerate}
\item[(a)] With probability $1 - 2\exp{(-\varepsilon^2/2)}$, $\norm{A_n^{-1}}_2 \leq 2\sqrt{n} + \varepsilon$,
\item[(b)] $\textnormal{E}(\norm{B_n^\textnormal{T} (A_n^{-1} b_n \odot A_n^{-\textnormal{T}} c_n)})_2^2 \geq \Omega(n)$,
\item[(c)] $\textnormal{E}(|c_n^\textnormal{T}A_n^{-1}b_n|^2) \leq O(1)$,
\end{enumerate}
and consequently, this family of problems on an average satisfies the $(1/2, 1/2)-$asymptotic convergence conditions.
\end{proposition}
Proposition 1 thus shows that theorem 1 is applicable on average for a design problem drawn from the distribution of problems defined in problem \ref{prob:gaussian_random}, and thus guarantees an average case convergence of gradient descent.

The second family of problems is one in which I fix the scaling of the singular values of the physics matrix with the state size, and generate the left and right singular vectors randomly --- the physical motivation behind this construction is that the scaling of the spectrum of the physics matrix with the problem size is often fixed by the derivative operators and boundary conditions appearing in the physical laws, while the precise singular vectors depend on the details of the (initial and randomly chosen) design parameters. Under some assumptions on the scalings of the singular values with the problem size, the asymptotic convergence conditions (definition~\ref{def:asym_conv_cond}) are shown to be satisfied on an average.

\begin{problem}
\label{prob:random_sing_vecs}
Given a cost function $f:\mathbb{R}\to \mathbb{R}$ being an $L-$Lipschitz smooth and $\mu-$strongly convex function that is bounded from below and $\gamma \geq 1/2$, define a family of physical design problems $\{(f, c_n, \varphi_n \equiv (A_n := R_n \textnormal{diag}(s_n)Q_n^\textnormal{T}, B_n:= I_n, b_n)) : n \in \mathbb{N}\}$ where
\begin{itemize}
\item $\norm{b_n}_\infty, \norm{b_n}_1, \norm{c_n}_\infty, \norm{c_n}_1 = \Theta(1)$,
\item $s_n \in (0, \infty)^n$ with $\norm{1 / s_n}_\infty \leq O(n^\gamma)$ and $\norm{1 / s_n} = \Theta(n)$,
\item  $R_n, Q_n \in \mathbb{R}^{n\times n}$ are drawn uniformly at random from the Haar measure over orthogonal matrices.
\end{itemize}
\end{problem}

\begin{proposition}  For the problem of state size $n$, $(f, c_n, \varphi_n \equiv (A_n , B_n, b_n))$, picked from the family of problems defined in problem \ref{prob:random_sing_vecs},
\begin{enumerate}
\item[(a)] $\norm{A_n^{-1}}_2 \leq  O(n^\gamma )$,
\item[(b)] $\textnormal{E}(\norm{B_n^\textnormal{T} (A_n^{-1} b_n \odot A_n^{-\textnormal{T}} c_n)})_2^2 \geq \Omega(n)$,
\item[(c)] $\textnormal{E}(|c_n^\textnormal{T}A_n^{-1}b_n|^2) \leq O(1)$,
\end{enumerate}
and consequently, this family of problems on an average satisfies the $1 -\gamma, \gamma$-asymptotic convergence conditions.
\end{proposition}
Consequently using theorem 1, it follows that for $\gamma \in [1/2, 3/4)$, a design problem drawn from the distribution of design problems defined in problem \ref{prob:gaussian_random} on an average is efficiently solvable by gradient descent.

Finally, I show that under some further assumptions on the physics matrix of the physical system, it can shown that all the designs generated during the gradient descent algorithm are non-resonant, and consequently gradient descent can be shown to converge to the global optima without the any additional assumptions on the gradient descent trajectory.

\begin{theorem}
Let $\mathcal{F} := \{(f, c_n, \varphi_n \equiv(A_n, B_n, b_n)) : n \in \mathbb{N}\}$ be a family of physical design problems that satisfies the $\alpha, \gamma-$asymptotic convergence conditions (definition 5) with $\alpha > 1/2, \gamma < 3\alpha$, and also satisfies $\norm{A_n^{-1}}_\textnormal{max} \leq O(n^{o(1)})$, then gradient descent when applied on $(f, c_n, \varphi_n)$, if it does not fail, produces a design $\theta_n^*$ such that $f(c^\textnormal{T} \varphi_n(\theta_n^*)) - f^* \leq \varepsilon$ in $T = \Theta(n^{1-2(\alpha - \gamma)}\log(\varepsilon^{-1}))$ steps.
\end{theorem}

The remainder of this paper is devoted to detailed proofs of these statements --- the proof of theorem 1 is provided in section \ref{sec:thm1}, section \ref{sec:prop12} is dedicated to the proofs of propositions 1 and 2 and the proof of theorem 2 is provided in section \ref{sec:thm2}.
\section{Detailed Proofs}
\subsection{Convergence of gradient descent}
\label{sec:thm1}
I begin by establishing an asymptotic property concerning the stability of a physical system --- the object of interest here is to study how a physical system behaves when the design parameters are perturbed slightly.
\begin{lemma}\label{lemma:perturb_norm}
Let $\{\varphi_n \equiv (A_n, B_n, b_n ): n \in \mathbb{N}\}$ be a sequence of physical systems which satisfies $\norm{B_n}_\infty \leq O(1)$ and $\norm{A^{-1}}_2 \leq O(n^\gamma)$ for some $\gamma > 0$, then for all sequences $\{\theta_n \in \textnormal{dom}(\varphi_n) : n \in \mathbb{N}\}$ such that $\norm{\theta_n}_\infty \leq O(n^{o(1) - \alpha - \gamma})$ for some $\alpha > 0$, $\norm{(A_n + \textnormal{diag}(B_n\theta_n))^{-1}}_2 \leq O(n^{\gamma})$.
\end{lemma}
\emph{Proof}: It follows straightforwardly that for $n \in \mathbb{N}$,
\begin{align*}
&(A_n + \textnormal{diag}(B_n \theta_n))^{-1} =\\
 &\qquad A_n^{-1} - A_n^{-1} \textnormal{diag}(B_n \theta_n) (A_n + \textnormal{diag}(B_n \theta_n))^{-1}.
\end{align*}
From the triangle inequality, I obtain that
\begin{align*}
&\norm{(A_n + \textnormal{diag}(B_n \theta_n))^{-1}}_2 \leq \norm{A_n^{-1}}_2 + \\
 &\qquad   \norm{A_n^{-1}}_2 \norm{B_n \theta_n}_\infty \norm{(A_n + \textnormal{diag}(B_n \theta_n))^{-1} }_2.
\end{align*}
By assumption, $\norm{A_n^{-1}}_2 \leq O(n^\gamma)$ for $\gamma > 0$, and $\norm{B_n \theta_n}_\infty \leq \norm{B_n}_\infty \norm{\theta_n}_\infty \leq O(n^{o(1) - \gamma - \alpha})$. Therefore,
\[
\big(1 - O(n^{o(1) - \alpha})\big) \norm{(A_n + \textnormal{diag}(B_n \theta_n))^{-1}}_2 \leq O(n^{\gamma}),
\]
from which the lemma statement follows. $\square$ 

For completeness, I provide an additional lemma with some standard properties of Lipschitz continuous and strongly convex functions that I will use in the following proofs.

\begin{lemma}
\label{lemma:fun_prop}
Let $f:\mathbb{R} \to \mathbb{R}$ be a twice-differentiable, $L-$Lipschitz smooth and $\mu-$strongly convex function such that $f^* = \min_{x\in \mathbb{R}}f(x) \geq -\infty$, then
\begin{enumerate}
\item[(a)] $\forall x, y \in \mathbb{R}$,
\[
f(x) \leq f(y) + \nabla f(y) (x - y) + \frac{L}{2} |x - y|^2.
\]
\item [(b)] $\forall x \in \mathbb{R}$,
\[
2\mu\big(f(x) - f^*\big) \leq  |\nabla f(x)|^2.
\]
\item[(c)] $\forall x \in \mathbb{R}$,
\[
|f'(x)|^2 \leq \frac{2L^2}{\mu} (f(x) - f^*).
\]
\end{enumerate}
\end{lemma}
\emph{Proof}:
\begin{enumerate}
\item[(a)] Since a $L-$Lipschitz smooth function by definition satisfies $| f'(x) - f'(y)| \leq L|x - y| \ \forall x, y \in \mathbb{R}$, it follows that $\forall x \in \mathbb{R}, \ | f''(x)| \leq L$. From Taylor's theorem, it follows that
\begin{align*}
&f(x) =\\
 &\quad f(y) + f'(y)(x - y) + \frac{1}{2}\int_{y}^x f''(s) (x - s) ds.
\end{align*}
Since $\forall s \in [x, y], \ f''(s) \leq |f''(s)| \leq L$, it follows that
\begin{align*}
&f(x) \leq f(y) + f'(y)(x - y) + \frac{L}{2}\int_{y}^x (x - s) ds \nonumber\\
&\qquad =f(y) + f'(y)(x - y) + \frac{L}{2}|x - y|^2.
\end{align*}
\item[(b)] Since $f$ is $\mu-$strongly convex, $\forall x, y \in \mathbb{R}$
\[
f(x) \geq f(y) + f'(y) (x - y) + \frac{\mu}{2}|x - y|^2.
\]
and hence $\forall y \in \mathbb{R}$
\[
\min_{x\in \mathbb{R}} f(x) \geq \min_{x \in \mathbb{R}} \bigg(f(y) + f'(y)(x - y) + \frac{\mu}{2}|x - y|^2\bigg),
\]
from which it follows that $\forall y \in \mathbb{R}$
\[
f^* \geq f(y) -\frac{1}{2\mu} | f'(y)|^2.
\]
\item[(c)] Let $x^* = \text{argmin}_{x\in \mathbb{R}}f^*$. I note that from stationarity conditions, $f'(x^*) = 0$. From the $L-$Lipschitz continuity, it follows that $|f'(x)|^2 \leq L^2(x - x^*)^2 \ \forall \ x \in \mathbb{R}$. Furthermore, from strong convexity of $f$, it follows that $\forall x \in \mathbb{R}$,
\[
f(x) \geq f^* + \frac{\mu}{2}|x - x^*|^2.
\]
Therefore, $\forall \ x\in\mathbb{R}$, $|f'(x)|^2 \leq 2L^2 / \mu (f(x) - f^*)$ $\square$.
\end{enumerate}

 Next, I analyze the gradient descent algorithm (algorithm \ref{alg:gd}). The next three lemmas characterize the decrease in the cost function on taking a gradient descent step.
\begin{lemma} 
\label{lemma:decomp}
Let $ (f, c, \varphi \equiv (A, B, b))$ be a physical design problem, and let $f$ be $L-$Lipschitz smooth, then $\forall \theta \in \textnormal{dom}(\varphi)$ and $\eta >0$ such that $\theta' := \theta - \eta \nabla_\theta f(c^\textnormal{T}\varphi(\theta)) \in \textnormal{dom}(\varphi)$,
\begin{align*}
&f(c^\textnormal{T}\varphi(\theta')) \leq  f(c^\textnormal{T}\varphi(\theta)) - \eta (f'(c^\textnormal{T}\varphi(\theta))^2\times \\
&\qquad \qquad\big(\norm{B^\textnormal{T} v(\theta)}_2^2 + \varepsilon_1(\theta, \theta') + \eta\varepsilon_2(\theta, \theta')\big),
\end{align*}
where $v(\theta) :=  \varphi(\theta)\odot \textnormal{ad}[\varphi](\theta, c)$
\begin{subequations}
\begin{align}
&\varepsilon_1(\theta, \theta') := \textnormal{ad}[\varphi](\theta, c)^\textnormal{T} \textnormal{diag}(BB^\textnormal{T} v(\theta)) (\varphi(\theta') - \varphi(\theta)), \\
&\varepsilon_2(\theta, \theta') := -\frac{L}{2}\big(\textnormal{ad}[\varphi](\theta, c)^\textnormal{T} \textnormal{diag}(BB^\textnormal{T}v(\theta)) \varphi(\theta')\big)^2.
\end{align}
\end{subequations}
\end{lemma}
\emph{Proof}: Since $f$ is $L-$Lipschitz smooth, it follows that
\begin{align}\label{eq:lipschitz_ineq}
&f(c^\textnormal{T}\varphi(\theta')) \leq f(c^\textnormal{T}\varphi(\theta)) + f'(c^\textnormal{T}\varphi(\theta)) c^\textnormal{T}(\varphi(\theta') - \varphi(\theta)) +\nonumber\\ &\qquad \qquad  \frac{L}{2}\big(c^\textnormal{T}(\varphi(\theta) - \varphi(\theta')\big)^2.
\end{align}
Furthermore, since
\[
\varphi(\theta') - \varphi(\theta) = -(A + \textnormal{diag}(B\theta))^{-1} \textnormal{diag}(B(\theta' - \theta)) \varphi(\theta'),
\]
I obtain that
\begin{align}\label{eq:olap_diff_1}
&c^\textnormal{T}(\varphi(\theta') - \varphi(\theta)) \nonumber \\
&=  -\textnormal{ad}[\varphi](\theta, c)^\textnormal{T} \varphi(\theta') \nonumber \\
&=-\eta f'(c^\textnormal{T}\varphi(\theta)) \textnormal{adj}[\varphi](\theta, c)^\textnormal{T}\textnormal{diag}(BB^\textnormal{T}v(\theta))\varphi(\theta').
\end{align}
Expressing $\varphi(\theta') = \varphi(\theta) + (\varphi(\theta') - \varphi(\theta))$, I obtain that
\begin{align}\label{eq:olap_diff_2}
&c^\textnormal{T}(\varphi(\theta') - \varphi(\theta)) = -\eta f'(c^\textnormal{T}\varphi(\theta))\big(\norm{B^\textnormal{T}v(\theta)}^2 + \nonumber\\
&\qquad \qquad \textnormal{adj}[\varphi](\theta, c)^\textnormal{T}\textnormal{diag}(BB^\textnormal{T}v(\theta))(\varphi(\theta') - \varphi(\theta))\big).
\end{align}
Substituting Eqs.~\ref{eq:olap_diff_1} and \ref{eq:olap_diff_2} into the Eq.~\ref{eq:lipschitz_ineq}, I obtain the lemma statement. $\square$. 
\begin{lemma}
\label{lemma:bound_eps_1}
Consider a family of physical design problems $\mathcal{F} := \{(f, c_n, \varphi_n \equiv (A_n, B_n, b_n)\}$ satisfying the $\alpha, \gamma$-asymptotic convergence conditions (definition \ref{def:asym_conv_cond}). Let $\eta > 0$, and let $\{\theta_n \in \textnormal{dom}(\varphi_n) : n \in \mathbb{N}\}$ be a sequence non-resonant with respect to $\mathcal{F}$ with $\norm{\theta_n}_\infty \leq O(n^{o(1) - \alpha - \gamma})$. Furthermore, suppose that the sequence $\{\theta_n' := \theta_n - \eta \nabla_{\theta}f(c^\text{T}_n \varphi_n(\theta_n)) : n \in \mathbb{N}\}$ satisfies $\norm{\theta_n'}_\infty \leq O(n^{o(1) - \alpha -\gamma})$, then $|\varepsilon_1(\theta_n, \theta_n')| \leq O(n^{o(1) + \gamma -\alpha + 1/2})$, where $\varepsilon_1$ is defined in lemma \ref{lemma:decomp}.
\end{lemma}
\emph{Proof}: From the definition of $\varepsilon_1$ it follows that $\forall n \in \mathbb{N}$,
\begin{align}\label{eq:epsilon_1_bound}
&|\varepsilon_1(\theta_n, \theta_n')| \leq \sqrt{n} \norm{\textnormal{ad}[\varphi_n](\theta_n, c_n)^\textnormal{T}\textnormal{diag}(B_nB_n^\textnormal{T}v_n(\theta_n))}_\infty \times \nonumber\\
&\qquad \qquad\norm{\varphi_n(\theta_n') - \varphi_n(\theta_n)}_2,
\end{align}
where $v_n(\theta_n) = \textnormal{ad}[\varphi_n](\theta_n, c_n) \odot \varphi_n(\theta_n)$. Since the sequence $\{\theta_n : n\in \mathbb{N}\}$ is non-resonant, $\norm{\textnormal{ad}[\varphi_n](\theta_n, c_n)^\textnormal{T}\textnormal{diag}(B_nB_n^\textnormal{T}v_n(\theta_n))}_\infty \leq O(n^{o(1)})$. Furthermore, $\forall n \in \mathbb{N}$
\begin{align*}
&\norm{\varphi_n(\theta_n) - \varphi_n(\theta_n')}_2 \leq \\
&\qquad\norm{(A_n + \textnormal{diag}(B_n \theta_n))^{-1}}_2 \norm{B_n(\theta_n - \theta_n')}_\infty \norm{\varphi_n(\theta_n')}_2
\end{align*}
It follows from lemma \ref{lemma:perturb_norm} that $\norm{(A_n + \textnormal{diag}(B_n\theta_n))^{-1}}_2 \leq O(n^\gamma)$ and $\norm{\varphi_n(\theta_n')}_2 \leq \norm{(A_n + \textnormal{diag}(B_n\theta_n))^{-1}}_2 \sqrt{\norm{b}_\infty \norm{b}_1} \leq O(n^{\gamma})$. Finally, $\norm{B_n(\theta_n - \theta_n')}_\infty \leq \norm{B}_n (\norm{\theta_n}_\infty + \norm{\theta_n'}_\infty) \leq O(n^{o(1) - \alpha -\gamma})$ and thus $\norm{\varphi_n(\theta_n') - \varphi_n(\theta_n)}_2 \leq O(n^{o(1) + \gamma - \alpha})$. Using these estimates the lemma statement follows.\\

\begin{lemma}
\label{lemma:bound_eps_2}
Consider a family of physical design problems $\mathcal{F}:=\{(f, c_n, \varphi_n \equiv (A_n, B_n, b_n))\}$ satisfying the $\alpha, \gamma-$asymptotic convergence conditions (definition \ref{def:asym_conv_cond}). Let $\eta > 0$ and $\{\theta_n \in \textnormal{dom}(\varphi_n) : n \in \mathbb{N}\}$ be a sequence non-resonant with respect to $\mathcal{F}$ with $\norm{\theta_n}_\infty \leq O(n^{o(1) - \alpha - \gamma})$. Furthermore, suppose that the sequence $\{\theta_n' := \theta_n - \eta \nabla_\theta f(c_n^\textnormal{T}\varphi_n(\theta_n)) : n \in \mathbb{N}\}$ satisfies $\norm{\theta_n' }_\infty \leq O(n^{o(1) - \alpha -\gamma})$, then $|\varepsilon_2(\theta_n, \theta_n')| \leq O(n^{o(1) + 4\gamma})$, where $\varepsilon_2$ is defined in lemma \ref{lemma:decomp}.
\end{lemma}
\emph{Proof}: I note that $\forall n \in \mathbb{N}$, $|\varepsilon_2(\theta_n)| \leq L/2 \norm{\textnormal{adj}[\varphi_n](\theta_n, c_n)} \norm{\varphi_n(\theta_n)}\norm{B_n B_n^\textnormal{T}v_n(\theta_n)}_\infty$, where, as in lemma \ref{lemma:bound_eps_1}, $v_n(\theta_n) := \textnormal{ad}[\varphi_n](\theta_n, c_n) \odot \varphi_n(\theta_n) \ \forall \ n \in \mathbb{N}$. It follows from lemma \ref{lemma:perturb_norm} and the observation that $\norm{b} \leq \sqrt{\norm{b}_1 \norm{b}_\infty} \leq O(1)$ that $\norm{\varphi_n(\theta_n)} \leq O(n^\gamma)$. Similarly, $\norm{\textnormal{ad}[\varphi_n](\theta_n, c_n)} \leq O(n^\gamma)$. Furthermore, $\norm{B_n B_n^\textnormal{T}}_\infty \leq \norm{B_n}_\infty \norm{B_n}_1 \leq O(1)$ and since $\{\theta_n: n \in \mathbb{N}\}$ is non-resonant, $\norm{\textnormal{ad}[\varphi_n](\theta_n, c_n)}_\infty \leq O(n^{o(1)})$. From these estimates, the lemma statement follows. \\

\begin{lemma}
\label{lemma:gradient_bound}
Consider a family of physical design problems $\mathcal{F} := \{(f, c_n, \varphi_n \equiv (A_n, B_n, b_n))\}$ satisfying the $\alpha, \gamma$-asymptotic convergence conditions (definition \ref{def:asym_conv_cond}) with $3\alpha > \gamma$ and $\norm{A_n^{-\textnormal{T}}c_n}_\infty \leq O(n^{o(1)})$. Let $\{\theta_n \in \textnormal{dom}(\varphi_n) : n\in \mathbb{N}\}$ be a sequence non-resonant with respect to $\mathcal{F}$ such that $\norm{\theta_n}_\infty \leq O(n^{o(1) - \alpha - \gamma})$, then $\norm{B_n^\textnormal{T} (\varphi_n(\theta_n) \odot \textnormal{ad}[\varphi_n](\theta_n, c_n))}^2 \geq \Omega(n^{\alpha + \gamma})$.
\end{lemma}
\emph{Proof}: For notational convenient, let $v_n(\theta) = \varphi_n(\theta) \odot \textnormal{ad}[\varphi_n](\theta, c_n) \ \forall \theta \in \textnormal{dom}(\varphi_n),  \ n \in \mathbb{N}$. I first upper bound $\norm{B_n^\textnormal{T}(v_n(\theta_n) - v_n(0))}$. Notice that
\[
\norm{B_n^\textnormal{T}(v_n(\theta_n) - v_n(0))} \leq \norm{B_n}_2 \norm{v_n(\theta_n) - v_n(0)}.
\]
Furthermore,
\begin{align*}
&\norm{v_n(\theta_n) - v_n(0)} \leq \nonumber \\
&\qquad \norm{\varphi_n(\theta_n)}_\infty \norm{\text{ad}[\varphi_n](\theta_n, c_n) - \text{ad}[\varphi_n](0, c_n)} + \\
&\qquad \norm{\textnormal{ad}[\varphi_n](0, c_n)}_\infty \norm{\varphi_n(\theta_n) - \varphi_n(0)}.
\end{align*}
Note that since $\{\theta_n: n\in \mathbb{N}\}$ is non-resonant, $\norm{\varphi_n(\theta_n)}_\infty \leq O(n^{o(1)})$ and by assumption $\norm{\textnormal{ad}[\varphi_n](0, c_n)}_\infty = \norm{A_n^{-\textnormal{T}}c_n}_\infty \leq O(n^{o(1)})$. Furthermore, $\forall n \in \mathbb{N}$
\[
\norm{\varphi_n(\theta_n) - \varphi_n(0)} \leq \norm{A_n^{-1}}_2 \norm{B_n\theta_n}_\infty \norm{\varphi_n(\theta_n)}.
\]
I note that by assumption $\norm{A_n^{-1}}_2 \leq O(n^\gamma)$, $\norm{B_n \theta_n} \leq O(n^{o(1) - \alpha - \gamma})$ and $\norm{\varphi_n(\theta_n)} \leq \norm{(A_n + \text{diag}(B_n \theta_n))^{-1}}_2 \norm{b_n} \leq O(n^\gamma)$. Consequently, $\norm{\varphi_n(\theta_n) - \varphi_n(0)} \leq O(n^{o(1) + \gamma - \alpha})$. Similarly, $\norm{\text{ad}[\varphi_n](\theta_n, c_n) - \text{ad}[\varphi_n](0, c_n)} \leq O(n^{o(1) + \gamma - \alpha})$, which yields
\[
\norm{B_n^\text{T}(v_n(\theta_n) - v_n(0))} \leq O(n^{o(1) + \gamma - \alpha}).
\]
Finally, I note that from the triangle inequality
\[
\norm{B_n^\text{T}v_n(\theta_n)} \geq \norm{B_n^\text{T}v_n(0)} - \norm{B_n^\text{T}(v_n(\theta_n) - v_n(0))}.
\]
Since, by assumption, $\norm{B_n^\text{T}v_n(0)} \geq \Omega(n^{(\alpha + \gamma)/2})$ and $\gamma < 3\alpha$, I obtain that $\norm{B_n^\text{T}v_n(\theta_n)} \geq \Omega(n^{(\alpha + \gamma)/2})$, from which the lemma statement follows. $\square$. \\

\noindent \textbf{Theorem 1 (Restated)} \emph{ Let $\mathcal{F} := \{(f, c_n, \varphi_n) : n \in \mathbb{N}\}$ be a family of physical design problems that satisfies the $\alpha, \gamma-$asymptotic convergence conditions (definition 5) with $\alpha > 1/4, \gamma < 3\alpha$, and for $n \in \mathbb{N}$, let $\Theta_n = \{\theta_{n}^1, \theta_{n}^2 \dots \theta_{n}^{T_n} \}$ be the design parameters generated by gradient descent (algorithm 1) under the assumption that it does not fail when applied on $(f, c_n, \varphi_n)$ with step size $\eta_n \leq O(n^{-3\gamma + \alpha - 1})$ for $T_n$ steps starting with initial design of $\theta_n^0 = 0$. If all the sequences $\{\theta_{n} \in \Theta_n: n \in \mathbb{N}\}$ are non-resonant with respect to $\mathcal{F}$, then $f(c_n^\textnormal{T}\varphi_n(\theta_{n}^{T_n})) - f^* \leq \varepsilon$ for $T$ chosen such that $T_n = \Theta(n^{1 - 2(\alpha - \gamma)}\log(\varepsilon^{-1}))$.} \\ \ \\
\noindent\emph{Proof}: Consider the trajectory $\{ \theta_{n}^1 \dots \theta_{n}^{T_n}\}$ generated by gradient descent when applied for $T_n$ steps on the problem $(f, c_n, \varphi_n)$ starting from $\theta_{n}^0 = 0$ and with step size $\eta_n$. I will first analyze the final cost function achieved under the assumption that $\forall t \in [T_n]$, $\norm{\theta_{n}^t}_\infty \leq O(n^{o(1) - \alpha -\gamma})$, and show that it can be made $\varepsilon-$close to $f^*$ after $T_n = \Theta(n^{1 -2(\alpha - \gamma)} \text{log}(\varepsilon^{-1}))$ gradient descent steps. Then, I will show that this assumption is valid for all gradient descent steps.

 I note with this assumption and the assumption that gradient descent algorithm generates a trajectory that is non-resonant, lemmas 4-7 are applicable. From lemmas \ref{lemma:bound_eps_1} and \ref{lemma:gradient_bound}, it is easy to see that if $\alpha > 1/4$, then $\forall t \in \{0, 1, 2 \dots T_n - 1\}$,
\[
\norm{B_n^\text{T}(\varphi_n(\theta_n^t)\cdot \text{ad}[\varphi_n](\theta_n^t, c_n)}^2 + 2 \varepsilon_1(\theta_n^t, \theta_n^{t + 1}) \geq \Omega(n^{\alpha + \gamma}).
\]
 Similarly, choosing $\eta_n \leq O(n^{-3\gamma + \alpha - 1})$, it follows from lemmas \ref{lemma:bound_eps_2} and \ref{lemma:gradient_bound} that $\forall t \in \{0, 1, 2 \dots T_n - 1\}$
 \[
 \norm{B_n^\text{T}(\varphi_n(\theta_n^t) \cdot \text{ad}[\varphi_n](\theta_n^t, c_n)}^2 + 2\varepsilon_2(\theta_n^t, \theta_n^{t + 1}) \geq \Omega(n^{\alpha + \gamma}).
 \]
 Consequently, from lemma \ref{lemma:decomp}, I obtain that $\forall t \in [T_n]$
\[
f(c_n^\text{T}\varphi_n(\theta_{t}^n)) \leq f(c_{n}^\text{T}\varphi_{n }(\theta^{t-1}_{n})) - \eta_nf'(c_n^\text{T}\varphi_n(\theta_n^{t-1}))^2  \Omega(n^{\alpha + \gamma}),
\]
From lemma \ref{lemma:fun_prop}b it follows that
\[
f(c_n^\text{T}\varphi_n(\theta_t^n)) - f^* \leq (f(c_n^\text{T}\varphi_n(\theta_{t-1}^n)) - f^*) (1 - \eta_n \Omega(n^{\alpha + \gamma})).
\]
where $f^* = \min_{x \in \mathbb{R}}f(x)$. Noting that, from lemma \ref{lemma:fun_prop}b along with asymptotic convergence condition (g), it follows that
\begin{align*}
f(c_n^\text{T}\varphi_n(\theta_0^n)) - f^* &= f(c_n^\text{T}A_n^{-1}b_n) - f^* \nonumber \\
&\leq \frac{L}{2}\big(c_n^\text{T}A_n^{-1}b_n - x^*)^2 \nonumber\\
&\leq O(1),
\end{align*}
where $x^* = \textnormal{argmin}_{x \in \mathbb{R}}f(x)$ and therefore
\[
f(c_n^\text{T}\varphi_n(\theta_{T_n}^n)) - f^* \leq O(1) (1 - \eta_n \Omega(n^{\alpha + \gamma}))^{T_n}.
\]
Consequently, using $T_n = \Theta(n^{1 - 2(\alpha - \gamma) }\text{log}(\varepsilon^{-1}))$, I obtain that $f(c_n^\text{T}\varphi_n(\theta^n_{T_n})) - f^* \leq \varepsilon$.

Next, I verify that the assumption $\norm{\theta_n^t}_\infty \leq O(n^{o(1) - \alpha - \gamma})$ holds for all $t \in [T_n]$. I notice that for all $t \in [T_n]$, the norm of $\theta_n^t$ can be upper bounded by the sum of norms of gradients in the previous steps, multiplied by the step size. Using lemma \ref{lemma:adj_method}
\begin{align*}
&\norm{\theta_n^t}_\infty \leq \nonumber\\
&\eta_n \sum_{t' = 0}^{t - 1}  |f'(c_n \varphi_n(\theta_n^{t'}))| \norm{B_n^\text{T}(\varphi_n(\theta_n^{t'}) \odot \text{ad}[\varphi_n](\theta_n^{t'}, c_n))}_\infty.
\end{align*}
Furthermore, using the fact that $\norm{B_n}_1 \leq O(1)$, and that the gradient descent trajectory is assumed to be non-resonant, $\norm{B_n^\text{T}(\varphi_n(\theta_n^{t}) \odot \text{ad}[\varphi_n](\theta_n^{t}, c_n))}_\infty \leq O(n^{o(1)})$ \ $\forall t \in [T_n]$. Furthermore, from lemma \ref{lemma:fun_prop}(c), it follows that
\begin{align*}
|f'(c_n \varphi_n(\theta_n^{t'}))| &\leq (2L^2 / \mu)^{1/2} (f(c_n^\text{T}\varphi_n(\theta_n^{t'})) - f^*))^{1/2} \nonumber\\
&\leq O(1)(1 - \eta_n \Omega(n^{\alpha +\gamma}))^{t'/2}.
\end{align*}
Therefore,
\begin{align*}
\norm{\theta_n^t}_\infty &\leq \eta_n O(n^{o(1)}) \sum_{t' = 0}^{t - 1} (1 - \eta_n \Omega(n^{\alpha +\gamma}))^{t'/2} \nonumber\\
&\leq \frac{\eta_n O(n^{o(1)})}{1 - (1 - \eta_n \Omega(n^{\alpha + \gamma}))^{1/2}}\nonumber\\
& \leq O(n^{o(1) - \alpha -\gamma}).
\end{align*}
This completes the proof of the lemma $\square$.

\subsection{Analysis of random family of physical design problems}
\label{sec:prop12}
\noindent\textbf{Problem 1 (Restated)} \emph{Given a cost function $f:\mathbb{R}\to \mathbb{R}$ being an $L-$Lipschitz smooth and $\mu-$strongly convex function that is bounded from below, define a family of physical design problems $\mathcal{F}:= \{(f, c_n, \varphi_n \equiv (A_n: = G_n^{-1},B_n:= I_n, b_n)) : n \in \mathbb{N}\}$ where
\begin{itemize}
\item $\norm{b_n}_\infty, \norm{b_n}_1, \norm{c_n}_\infty, \norm{c_n}_1 = \Theta(1)$ and,
\item $G_n \in \mathbb{R}^{n\times n}$ is a matrix where each entry is independently drawn from the standard normal distribution.
\end{itemize}} 
\noindent\textbf{Proposition 1 (Restated)} \emph{For the problem of state size $n$, $(f, c_n, \varphi_n \equiv (A_n, B_n, b_n))$, picked from the family of problems defined in problem \ref{prob:gaussian_random}, it is true that
\begin{enumerate}
\item[(a)] With probability $1 - 2\exp{(-\varepsilon^2/2)}$, $\norm{A_n^{-1}}_2 \leq 2\sqrt{n} + \varepsilon$,
\item[(b)] $\textnormal{E}(\norm{B_n^\textnormal{T} (A_n^{-1} b_n \odot A_n^{-\textnormal{T}} c_n)})_2^2 \geq \Omega(n)$,
\item[(c)] $\textnormal{E}(|c_n^\textnormal{T}A_n^{-1}b_n|^2) \leq O(1)$,
\end{enumerate}
and consequently, this family of problems on an average satisfies the $1/2,1/2-$asymptotic convergence conditions.} \\

\noindent\emph{Proof}: 
\begin{enumerate}
\item[(a)] This is a standard result in random matrix theory, for e.g.~see Ref.~\cite{davidson2001local}.
\item[(b)] It follows from straightforward computation that
\begin{align*}
&\textnormal{E}(\norm{B_n^\textnormal{T} (A_n^{-1} b_n \odot A_n^{-\textnormal{T}} c_n)})_2^2 \nonumber \\
&\qquad= \text{E}(\norm{G_n b_n \odot G_n^\text{T} c_n}^2) \\
&\qquad= n \norm{b_n}^2 \norm{c_n}^2 + \norm{b_n\odot c_n}^2.
\end{align*}
Noting that $\norm{b_n}^2 \geq \norm{b_n}_\infty^2 \geq \Omega(1)$ and $\norm{c_n}^2 \geq \norm{c_n}_\infty^2 \geq \Omega(1)$, the lemma statement follows.
\item[(c)] I note that. 
\begin{align*}
\text{E}(|c_n^\text{T}A_n^{-1}b_n|^2) &=  \text{E}\big(c_n^\text{T}G_n b_n b_n^\text{T}G_n^\text{T}c_n\big) \nonumber\\
&= \norm{c_n}^2\norm{b_n}^2.
\end{align*}
Since $\norm{b_n}^2 \leq \norm{b_n}_1^2 \leq O(1)$ and $\norm{c_n}^2 \leq \norm{c_n}_1^2 \leq O(1)$, the lemma statement follows. $\square$
\end{enumerate}

\noindent \textbf{Problem 2 (Restated)} \emph{Given a cost function $f:\mathbb{R}\to \mathbb{R}$ being an $L-$Lipschitz smooth and $\mu-$strongly convex function that is bounded from below and $\gamma \geq 1/2$, define a family of physical design problems $\{(f, c_n, \varphi_n \equiv (A_n := R_n \textnormal{diag}(s_n)Q_n^\textnormal{T}, B_n:= I_n, b_n)) : n \in \mathbb{N}\}$ where
\begin{itemize}
\item $\norm{b_n}_\infty, \norm{b_n}_1, \norm{c_n}_\infty, \norm{c_n}_1 = \Theta(1)$,
\item $s_n \in (0, \infty)^n$ with $\norm{1 / s_n}_\infty \leq O(n^\gamma)$ and $\norm{1 / s_n} = \Theta(n)$,
\item  $R_n, Q_n \in \mathbb{R}^{n\times n}$ are drawn uniformly at random from the Haar measure over orthogonal matrices.
\end{itemize}}
\begin{lemma}
\label{lemma:haar_second_moment}
Let $M = R\ \textnormal{diag}(v) Q^\text{T} \in \mathbb{R}^{d\times d}$, where $R, Q$ are orthogonal matrices drawn independently from the Haar random measure over the set of $d\times d$ orthogonal matrices, and $b, c \in \mathbb{R}^d$, then
\[
\textnormal{E}\big((c^\textnormal{T}M b)^2\big) = \frac{1}{d^2} \norm{v}^2 \norm{b}^2 \norm{c}^2.
\]
\end{lemma}
\emph{Proof}: Provided in appendix \ref{app:orthogonal_mean}.
\ \\
\begin{lemma}
\label{lemma:haar_fourth_moment}
Let $M = R\ \textnormal{diag}(v) Q^\text{T} \in \mathbb{R}^{d\times d}$, where $R, Q$ are orthogonal matrices drawn independently from the Haar random measure over the set of $d\times d$ orthogonal matrices, and $b, c \in \mathbb{R}^d$, then
\[
\textnormal{E}\big(\norm{Mb \odot M^\textnormal{T}c}^2\big) = \norm{b}^2 \norm{c}^2 \norm{v}^4 \Omega\bigg(\frac{1}{d^3}\bigg).
\]
\end{lemma}
\emph{Proof}: Provided in appendix \ref{app:orthogonal_mean}.
\ \\ \ \\
\noindent \textbf{Proposition 2 (Restated)} \emph{ For the problem of state size $n$, $(f, c_n, \varphi_n \equiv (A_n , B_n, b_n))$, picked from the family of problems defined in problem \ref{prob:random_sing_vecs},
\begin{enumerate}
\item[(a)] $\norm{A_n^{-1}}_2 \leq  O(n^\gamma)$,
\item[(b)] $\textnormal{E}(\norm{B_n^\textnormal{T} (A_n^{-1} b_n \odot A_n^{-\textnormal{T}} c_n)})^2 \geq \Omega(n)$,
\item[(c)] $\textnormal{E}(|c_n^\textnormal{T}A_n^{-1}b_n|^2) \leq O(1)$,
\end{enumerate}
and consequently, this family of problems on an average satisfies the $1 - \gamma, \gamma$-asymptotic convergence conditions.
}
\ \\ \ \\
\noindent\emph{Proof}:
\begin{itemize}
\item[(a)] This follows straightforwardly by noting that $\norm{R_n} = \norm{Q_n} = 1 \ \forall \ n \in \mathbb{N}$ and $\norm{\text{diag}(s_n)^{-1}} = \norm{1 / s_n}_\infty \leq O(n^\gamma)$.
\item[(b)] Using lemma \ref{lemma:haar_fourth_moment}, it immediately follows that
\begin{align*}
&\text{E}(\norm{B_n^\text{T}(A_n^{-1}b_n \odot A_n^{-\text{T}}c_n)}^2) \nonumber\\
&\qquad\geq \norm{b_n}^2\norm{c_n}^2 \norm{1/s_n}^4 \Omega(n^{-3}).
\end{align*}
Noting that $\norm{b_n}^2 \norm{c_n}^2 \geq \norm{b_n}_\infty^2\norm{c_n}_\infty^2 \geq \Omega(1)$ and $\norm{1 / s_n}^4 \geq \Omega(n^4)$, the lemma statement follows.
\item [(c)]  Using lemma \ref{lemma:haar_second_moment}, it follows that
\begin{align*}
\text{E}\big((c_n^\text{T}A_n^{-1}b_n)^2\big) = \frac{1}{n^2}\norm{c_n}^2\norm{b_n}^2\norm{1/s_n}^2.
\end{align*}
Since $\norm{c_n}^2\norm{b_n}^2 \leq \norm{c_n}_1^2 \norm{b_n}_1^2\leq O(1)$, and $\norm{1 / s_n}^2 \leq O(n^2)$, the lemma statement follows. $\square$
\end{itemize}
\subsection{Provably non-resonant problems}
\label{sec:thm2}
\begin{lemma}\label{lemma:perturb_max}
Let $\{\varphi_n \equiv (A_n, B_n, b_n): n \in \mathbb{N}\}$ be a sequence of physical systems which satisfies $\norm{b_n}_1 \leq O(1)$, $\norm{B_n}_\infty \leq O(1)$, $\norm{A_n^{-1}}_\textnormal{max} \leq O(n^{o(1)})$ and $\norm{A^{-1}}_2 \leq O(n^\gamma)$ for some $\gamma < 1/2$, then $\forall$ sequences $\{\theta_n \in \textnormal{dom}(\varphi_n) : n \in \mathbb{N}\}, \{c_n \in \mathbb{R}^n : n\in \mathbb{N}\}$ such that $\norm{\theta_n}_\infty \leq O(n^{o(1) - 1})$ and $\norm{c_n}_1 \leq O(1)$, $\norm{\varphi_n(\theta_n)}_\infty \leq O(n^{o(1)})$ and $\norm{\textnormal{ad}[\varphi_n](\theta_n, c_n)}_\infty \leq O(n^{o(1)})$.
\end{lemma}
\noindent\emph{Proof}: It follows straightforwardly that for $n \in \mathbb{N}$,
\begin{align*}
\varphi_n(\theta_n) = A_n^{-1}b_n - A_n^{-1} \textnormal{diag}(B_n \theta_n) \varphi_n(\theta_n).
\end{align*}
Consequently,
\begin{align*}
&\norm{\varphi_n(\theta_n)}_\infty \leq \norm{A_n^{-1}b_n}_\infty + \norm{A_n^{-1} \textnormal{diag}(B_n \theta_n) \varphi_n(\theta_n)}_\infty\\
 &\qquad \leq  \norm{A_n^{-1}}_\textnormal{max}\norm{b_n}_1 + \sqrt{n} \norm{A_n^{-1}}_2 \norm{B_n \theta_n}_\infty \norm{\varphi_n(\theta_n)}_\infty.
\end{align*}
By assumption, $\norm{A_n^{-1}}_2 \leq O(n^\gamma)$ for $\gamma < 1/2$, and $\norm{B_n \theta_n}_\infty \leq \norm{B_n}_\infty \norm{\theta_n}_\infty \leq O(n^{o(1) - 1})$. Therefore,
\[
\big(1 - O(n^{o(1) + \gamma - 1/2})\big) \norm{\varphi_n(\theta_n)}_\infty \leq O(n^{o(1)}),
\]
from which it follows that $\norm{\varphi_n(\theta_n)}\leq O(n^{o(1)})$. A similar analysis yields $\norm{\textnormal{ad}[\varphi_n](\theta_n, c_n)}_\infty \leq O(n^{o(1)})$. $\square$ \\

\noindent\textbf{Theorem 2 (Restated)} \emph{Let $\mathcal{F} := \{(f, c_n, \varphi_n \equiv(A_n, B_n, b_n)) : n \in \mathbb{N}\}$ be a family of physical design problems that satisfies the $\alpha, \gamma-$asymptotic convergence conditions (definition 5) with $\alpha > 1/2, \gamma < 3\alpha$, and also satisfies $\norm{A_n^{-1}}_\textnormal{max} \leq O(n^{o(1)})$, then gradient descent when applied on $(f, c_n, \varphi_n)$, if it does not fail, produces a design $\theta_n^*$ such that $f(c^\textnormal{T} \varphi_n(\theta_n^*)) - f^* \leq \varepsilon$ in $T = \Theta(n^{1-2(\alpha - \gamma)}\log(\varepsilon^{-1}))$ steps.} \ \\

\noindent\emph{Proof}: This theorem can be proved by repeating the analysis in the proof of theorem 1, and noting from lemma \ref{lemma:perturb_max} that the assumption that the gradient descent trajectory $\{\theta_n^1, \theta_n^2 \dots \}$ obtained for the problem $(f, c_n, \varphi_n)$ satisfies $\norm{\theta_n^t}_\infty \leq O(n^{\alpha + \gamma})$ implies that the trajectory is also non-resonant. $\square$
\section{Conclusion and open problems}
In conclusion, this work provides rigorous evidence for local optimization algorithms being efficient at solving physical design problems. I show that, under some assumptions on the physics of the system, non-resonant physical design problems are efficiently solvable by gradient descent. Furthermore, I also outline random ensembles of physical design problems which are, on an average, efficiently solvable by local optimization algorithms.

This work, while being a first step towards theoretically understanding the complexity of typical physical design problems, leaves several questions open. One question is to better characterize when a physical design problem is resonant --- in this analysis, I need to assume that gradient descent avoids resonant designs in order to show that it globally solves the design problem. While theorem 2 makes some progress in this direction, it would be interesting to more carefully analyze the gradient descent trajectory and obtain a set of weaker conditions on the design problems under which gradient descent avoids resonant devices. Another interesting direction would be to study the optimality of gradient descent on specific design problems, or ensembles of design problems, appearing in practical settings using the tools introduced in this paper. Finally, extending the analysis introduced in this paper to other optimization algorithms (such as quasi-Newton methods like BFGS, L-BFGS, or method of moving asymptotes), would also go a long way in making the rigorous results practically relevant.

\begin{acknowledgements}
I thank Shivam Garg, Logan Su, Geunho Ahn and Alex White for useful discussion. I acknowledge support from Max Planck Harvard Research Center for Quantum Optics (MPHQ) postdoctoral fellowship.
\end{acknowledgements}
\bibliography{references.bib}
\onecolumngrid
\appendix
\section{Proof of lemmas}
\label{app:orthogonal_mean}
This proof requires a basic integration formula with respect to the Haar measure over the orthogonal matrices. Given that $O \in \mathbb{R}^{n\times n}$ is an orthogonal matrix drawn randomly from the Haar measure, this integration formula allows us to express $\text{E}(O_{i_1, j_1} O_{i_2, j_2} \dots O_{i_{2k}, j_{2k}})$ in terms of the orthogonal Wein-garten function \cite{collins2006integration, collins2009some}. This function is difficult to evaluate in closed form for a general $k$, but for my purposes $k = 2$ will suffice. For completeness, I provide this integration formula, and specialize it to $k = 2$ and then provide a proof of proposition 2(b).
\begin{definition}[Pairing] For $n \in \mathbb{N}$, denoted by $\mathcal{P}_{2n}$, is an unordered tuple of $n$ unordered tuples with two elements, $((i_1, j_1), (i_2, j_2) \dots (i_n, j_n))$ where all $i_1, i_2 \dots i_n, j_1, j_2 \dots j_n$ are distinct and $\in [2n]$. The set of all pairings over $[2n]$ will be denoted by $\mathcal{P}_{2n}$.
\end{definition}
\textbf{Remarks}:
\begin{itemize}
\item It is important to emphasize the simply reordering the tuples in a pairing, or the the two elements inside the tuple, does not generate a different pairing. For instance, the pairing $((1, 2), (3, 4), (5, 6))$ over $[6]$ is the same as the pairing $((3, 4), (2, 1), (6, 5))$.
\item As an explicit and important example for the following calculations, the set $\mathcal{P}_4$ has three distinct elements, $((1, 2), (3, 4))$, $((1, 3), (2, 4))$, $((1, 4), (2, 3))$.
\end{itemize}
\begin{definition}[Pairing delta function $\Delta^p_{i_1, i_2 \dots i_{2n}}$]
\label{def:pairing_delta_fun}
Given a pairing $p = ((k_1, r_1), (k_2, r_2) \dots (k_n, r_n)) \in \mathcal{P}_{2n}$, and indices $i_1, i_2 \dots i_{2n}$ from some set $\mathcal{I}$, then
\[
\Delta^{p}_{i_1, i_2 \dots i_{2n}} = \prod_{s = 1}^{n} \delta_{i_{k_s}, i_{r_s}},
\]
where $\delta_{i, j} = 1$ if $i = j$ and $0$ if $i \neq j$ is the Kronecker delta function.
\end{definition}
\begin{definition}[Loop function $\textnormal{loop}(p_1, p_2)$] Let $p_1, p_2 \in \mathcal{P}_{2k}$ be two pairings. Construct a graph with vertices $\{1, 2, 3 \dots 2k\}$ with edges on all the pairings in $p_1$ or $p_2$, then $\textnormal{loop}(p_1, p_2)$ is the number of connected components in the graph.
\end{definition}
\begin{definition}[Orthogonal wein-garten function $\textnormal{Wg}_k^{O(d)}(p_1, p_2)$] Given $\mathcal{P}_{2k} = \{p_1, p_2 \dots\}$, let $G_{k}^d \in \mathbb{R}^{|\mathcal{P}_{2k}| \times |\mathcal{P}_{2k}|}$ be a matrix with elements $(G_{k}^d)_{i, j} = d^{\textnormal{loop}(p_i, p_j)}$, then the $\textnormal{Wg}_k^{O(d)}(p_i, p_j) = [(G_{k}^d)^{-1}]_{i, j}$.
\end{definition}
\begin{lemma}
\label{lemma:wg_4}
For $d > 1$ and for $p_1, p_2 \in \mathcal{P}_{4}$,
\begin{align*}
\textnormal{Wg}_2^{O(d)}(p_1, p_2) = 
\begin{cases}
\frac{1}{d^2 - 1} & \text{ if } p_1 = p_2, \\
-\frac{1}{d(d^2 - 1)} & \text{ if } p_1 \neq p_2.
\end{cases}
\end{align*}
\end{lemma}
\emph{Proof}: Explicitly, $\mathcal{P}_4 = \{p_1 := ((1, 2), (3, 4)), p_2:= ((1, 3), (2, 4)), p_3:= ((1, 4), (2, 3))\}$. I can now calculate the loop function, and it follows that $\text{loop}(p_i, p_j) = 2$ if $i = j$ else $1$ and thus
\[
G_2^d = \begin{bmatrix} d^2 & d & d \\
d & d^2 & d \\
d & d & d^2
\end{bmatrix} \implies
\big(G_2^d \big)^{-1} =\frac{1}{d^2 - 1} \begin{bmatrix} 1 & -1/d & -1/d \\
-1/d & 1 & -1/d \\
-1/d & -1/d & 1
\end{bmatrix}.
\] 
Identifying the weingarten function with these matrix elements, the lemma statement follows $\square$.

\begin{lemma}[Integration w.r.t. Haar measure over $O(d)$ from Ref.~\cite{collins2009some}]
\label{lemma:int_formula}
 Let $R$ be a matrix drawn uniformly at random from the Haar measure over the set of $d\times d$ orthogonal matrices, then for all $k \in \mathbb{N}$,
\[
\text{E}\big[R_{i_1, j_1}R_{i_2, j_2}, R_{i_3, j_3} \dots R_{i_{2k}, j_{2k}}\big] = \sum_{p_1, p_2 \in \mathcal{P}_{2k}} \textnormal{Wg}^{O(d)}_k(p_1, p_2) \Delta^{p_1}_{i_1, i_2 \dots i_{2k}} \Delta^{p_2}_{j_1, j_2 \dots j_{2k}}.
\]
\end{lemma}
\textbf{Lemma 8 (Restated)} \emph{Let $M = R\ \textnormal{diag}(v) Q^\text{T} \in \mathbb{R}^{d\times d}$, where $R, Q$ are orthogonal matrices drawn independently from the Haar random measure over the set of $d\times d$ orthogonal matrices, and $b, c \in \mathbb{R}^d$, then
\[
\textnormal{E}\big((c^\textnormal{T}M b)^2\big) = \frac{1}{d^2} \norm{v}^2 \norm{b}^2 \norm{c}^2.
\]
}
\emph{Proof}: Explicitly writing out the expectation value,
\[
\textnormal{E}((c^\text{T}Mb)^2) = \textnormal{E}\bigg(\sum_{i, j, k \in [d]^2}\prod_{l=1}^{2}v_{k_l}  c_{i_l}  b_{j_l}  R_{i_l, k_l} Q_{j_l, k_2}\bigg) = \sum_{i, j, k \in [d]^2}\bigg(\prod_{l=1}^{2}v_{k_l}  c_{i_l}  b_{j_l}\bigg)   \textnormal{E}\big(R_{i_1, k_1} R_{i_2, k_2}\big) \text{E}\big(Q_{j_1, k_1} Q_{j_2, k_2}\big).
\]
From lemma \ref{lemma:int_formula}, it follows that
\[
\textnormal{E}(R_{i_1, k_1} R_{i_2, k_2}) = \frac{1}{d} \delta_{i_1, i_2}\delta_{k_1, k_2}, \text{ and } \textnormal{E}(Q_{j_1, k_1} Q_{j_2, k_2}) = \frac{1}{d} \delta_{j_1, j_2}\delta_{k_1, k_2},
\]
and consequently,
\[
\text{E}\big((c^\text{T}Mb)^2\big) = \frac{1}{d^2}\norm{v}^2\norm{b}^2\norm{c}^2. 
\]
and thus the lemma statement follows $\square$. \\ \ \\
\textbf{Lemma 9 (Restated)}\emph{ Let $M = R\ \textnormal{diag}(v) Q^\text{T} \in \mathbb{R}^{d\times d}$, where $R, Q$ are orthogonal matrices drawn independently from the Haar random measure over the set of $d\times d$ orthogonal matrices, and $b, c \in \mathbb{R}^d$, then
\[
\textnormal{E}\big(\norm{Mb \odot M^\textnormal{T}c}^2\big) = \norm{b}^2 \norm{c}^2 \norm{v}^4 \Omega(d^{-3}).
\]
}
\emph{Proof}: It is easily follows that
\begin{align*}
\textnormal{E}\big(\norm{Mb \odot M^\textnormal{T}c}^2\big) = \sum_{i \in [n]} \sum_{j, k, l, m \in [n]^2} v_{j_1}v_{j_2} v_{k_1}v_{k_2} \text{E}\big(R_{i, j_1} R_{i, j_2} R_{l_1, k_1} R_{l_2, k_2}  \big) \text{E}\big(Q_{m_1, j_1} Q_{m_2, j_2}Q_{i, k_1} Q_{i, k_2}  \big) b_{m_1}b_{m_2} c_{l_1} c_{l_2}.
\end{align*}
From lemma \ref{lemma:int_formula}, it follows that
\begin{align*}
&\text{E}\big(R_{i, j_1} R_{i, j_2} R_{l_1, k_1} R_{l_2, k_2}  \big) = \sum_{p_1, p_2 \in \mathcal{P}_4} \text{Wg}^{O(d)}_4(p_1, p_2) \Delta^{p_1}_{i, i, l_1, l_2} \Delta^{p_2}_{j_1, j_2, k_1, k_2}, \\
&\text{E}\big(Q_{m_1, j_1} Q_{m_2, j_2}Q_{i, k_1} Q_{i, k_2} \big) = \sum_{p_1', p_2' \in \mathcal{P}_4} \text{Wg}^{O(d)}_4(p_1', p_2') \Delta^{p_1'}_{m_1, m_2, i, i} \Delta^{p_2'}_{ k_1, k_2, j_1, j_2}.
\end{align*}
I first evaluate the summation
\begin{align*}
f(p_1, p_1', v) := \sum_{j, k \in [d]^2}\sum_{p_2, p_2' \in \mathcal{P}_4} v_{j_1}v_{j_2} v_{k_1}v_{k_2} \text{Wg}_4^{O(d)}(p_1, p_2) \text{Wg}_4^{O(d)}(p_1', p_2')  \Delta^{p_2}_{j_1, j_2, k_1, k_2}\Delta^{p_2'}_{j_1, j_2, k_1, k_2}.
\end{align*}
I note from definition \ref{def:pairing_delta_fun}
\[
\sum_{j, k \in [d]^2}v_{j_1}v_{j_2} v_{k_1}v_{k_2} \Delta^{p_2}_{j_1, j_2, k_1, k_2} \Delta^{p_2'}_{k_1, k_2, j_1, j_2} = \begin{cases}
\norm{v}_4^4 & \text{ if } p_2 \neq p_2', \\
\norm{v}^4 & \text{ if } p_2 = p_2'.
\end{cases}
\]
Consequently, it follows that
\begin{align*}
f(p_1, p_1', v) = \big(\norm{v}^4 - \norm{v}_4^4\big) \sum_{p\in \mathcal{P}_4}\text{Wg}_4^{O(d)}(p_1, p) \text{Wg}_4^{O(d)}(p_1', p) + \norm{v}_4^4 \bigg(\sum_{p \in \mathcal{P}_4} \text{Wg}_4^{O(d)}(p_1, p)\bigg) \bigg(\sum_{p \in \mathcal{P}_4} \text{Wg}_4^{O(d)}(p_1', p)\bigg).
\end{align*}
Using the explicit formula for the Weingarten functions from from lemma \ref{lemma:wg_4}, it follows that
\[
\chi(p_1, p_1', v) = \begin{cases} \chi_\text{eq}(v)  & \text{ if } p_1 = p_1', \\
\chi_\text{ueq}(v) & \text{ if } p_1\neq p_1'.
\end{cases}, 
\]
where
\begin{align*}
& \chi_\text{eq}(v) = \frac{d^2 + 2}{d^2(d^2 - 1)^2} \norm{v}^4 - \frac{4d - 2}{d^2(d^2 - 1)^2} \norm{v}_4^4,\\
& \chi_\text{ueq}(v) = -\frac{2d - 1}{d^2(d^2 - 1)^2} \norm{v}^4 + \frac{d^2 - 2d + 3}{d^2(d^2 - 1)^2} \norm{v}_4^4
\end{align*}
I note that
\begin{align*}
\text{E}\big(\norm{Mb \odot M^\text{T}c}^2\big) &= \sum_{p_1, p_1' \in \mathcal{P}_4} \sum_{\substack{i \in [n],\\ l, m \in [n]^2}} b_{m_1} b_{m_2} c_{l_1} c_{l_2} \Delta^{p_1}_{i, i, l_1, l_2}\Delta^{p_1'}_{m_1, m_2, i, i} \chi(p_1, p_1', v) \nonumber\\
&=\chi_\text{eq}(v) \big(d\norm{b}^2 \norm{c}^2 + 2\norm{b\odot c}^2\big) + 6\chi_\text{ueq}(v) \norm{b}^2 \norm{c}^2.
\end{align*}
Using the expressions for $\chi_\text{eq}(v), \chi_\text{uneq}(v)$, we prove the lemma statement. $\square$.

\end{document}